\documentclass[11pt]{article}
\usepackage[latin1]{inputenc}
\usepackage{amsfonts}
\usepackage{amssymb}
\newcommand{\be}{\begin{equation}}
\newcommand{\ee}{\end{equation}}
\newcommand{\ba}{\begin{array}}
\newcommand{\ea}{\end{array}}
\newcommand{\bea}{\begin{eqnarray}}
\newcommand{\eea}{\end{eqnarray}}
\newcommand{\bee}{\begin{eqnarray*}}
\newcommand{\eee}{\end{eqnarray*}}

\newtheorem{Rk}{Remark}

\newenvironment{proof}{\noindent{\bf Proof.~}}
{{\mbox{}\hfill {\small \fbox{}}\\}}

\catcode`@=11
\renewcommand\appendix{\bigskip {\noindent\Large \bf Appendix}\par
  \setcounter{section}{0}%
  \setcounter{subsection}{0}%
  \renewcommand\thesection{\@Alph\c@section}}
\catcode`@=12

\def\Section{\setcounter{equation}{0}\section}
\newtheorem{theorem}{Theorem}[section]

\newtheorem{proposition}[theorem]{Proposition}

\def\thesection{\arabic{section}}

\newenvironment{acknowledgement}{\noindent{\bf Acknowledgement.~}}{}

\def\RR{\mathbb{R}}
\def\R{{\mathbb R}}

\def\lim{\mathop{\rm lim}}

\def\sup{\mathop{\rm sup}}

\def\e{\varepsilon}

\def\fref#1{{\rm (\ref{#1})}}

\def\ds{\displaystyle}
\def\ni{\noindent}
\def\bs{\bigskip}

\def\calA{{\cal A}}

\def\calE{{\cal E}}
\def\calH{{\cal H}}
\def\calC{{\cal C}}
\def\calF{{\cal F}}

\def\calS{{\cal S}}

\def\eps{\varepsilon}
\def\alphadef{m}
\def\ds{\displaystyle}
\def\pa{\partial}
\def\qed{\mbox{}\hfill {\small \fbox{}}\\}

\title{Stable ground states for the relativistic gravitational Vlasov-Poisson system}
\date{}

\author{Mohammed Lemou$^{**,***}$, Florian M\'ehats$^{**}$ and Pierre
  Rapha\"el$^{*}$}

\date{\it $^{*}$ IMT, Universit\'e Paul Sabatier, Toulouse, France\\ 
$^{**}$ IRMAR, Universit\'e Rennes 1, France\\
$^{***}$ CNRS, France.}

\begin{document}

\maketitle

\centerline{\sl This work was submitted on January 31, 2008.}

\bs


\begin{abstract}
We consider the three dimensional gravitational Vlasov-Poisson (GVP) system
in both classical and relativistic cases. The classical problem is subcritical in the natural energy space and the stability of a large class of ground states has been derived by various authors. The relativistic problem is critical and displays finite time blow up solutions. Using standard concentration compactness techniques, we however show that the breaking of the scaling symmetry allows the existence of stable relativistic ground states. A new feature in our analysis which applies both to the classical and relativistic problem is that the orbital stability of the ground states does not rely as usual on an argument of uniqueness of suitable minimizers --which is mostly unknown-- but on strong rigidity properties of the transport flow, and this extends the class of minimizers for which orbital stability is now proved.

\end{abstract}


%


\Section{Introduction}


\subsection{Setting of the problem}


We consider the three dimensional gravitational Vlasov-Poisson (GVP) system

\be
\label{rvp}
\left   \{ \begin{array}{llll}
         \ds \pa_t f+\frac{v}{\sqrt{1+|v|^2/c^2}}\cdot\nabla_x f-\nabla
         \phi_f \cdot\nabla_v f=0, \ \ (t,x,v)\in \RR_+\times \RR^3\times \RR^3\\[3mm]
          \ds f(t=0,x,v)=f_0(x,v)\geq 0,\\[3mm]
         \end{array}
\right.
\end{equation}
where
\be
\label{phif}
 \phi_f(x)=-\frac{1}{4\pi |x|}\star \rho_f,\ \ \rho_f(x)=\int_{\RR^3} f(x,v)\,dv,
\end{equation}
and $c\in ]0,+\infty]$ is the dimensionless light speed. The value $c=+\infty$
    recovers the classical Vlasov-Poisson system,
  which is a nonlinear transport equation describing  the mechanical state of a
stellar system subject to its own gravity (see for instance
\cite{binney,fridman}). In some situations when high velocities
occur, relativistic corrections should be introduced, see Van Kampen
and Felderhof \cite{VKF}, Glassey and Schaeffer \cite{GS} and
references therein. A more accurate model is then provided by
\fref{rvp} with $0<c<+\infty$, which is  the so called three dimensional relativistic gravitational Vlasov-Poisson system.\\

These systems are Hamiltonian and  all smooth enough solutions to
\fref{rvp} satisfy the  conservation of the $L^q$ norm
and of the total energy (Hamiltonian) on their lifespan:

\begin{equation}
\label{loi1}
\forall t,  \quad\forall q\in [1,+\infty], \ \ \quad |f(t)|_{L^q}=|f_0|_{L^q},\ \ {\calH_c}(f(t))={\calH_c}(f(0)),
\end{equation}
with
\begin{equation}
\label{defhamiltonian}
{\cal{H}}_c(f(t))=  
          \ds \int_{\RR^6}
         \gamma_c(v)  f(t,x,v)dxdv-\frac{1}{2}\int_{\RR^3}|\nabla
          \phi_f(t,x)|^2dx,
\end{equation}
\be
\label{gammac}
\ds \gamma_c(v)=    c^2\left(
          \sqrt{1+|v|^2/c^2} -1\right). 
\end{equation}
In particular
\begin{equation}
\label{defhamiltonian-infini}
{\cal{H}}_{\infty}(f(t))=  
          \ds \frac{1}{2}\int_{\RR^6}
         |v|^2  f(t,x,v)dxdv-\frac{1}{2}\int_{\RR^3}|\nabla
          \phi_f(t,x)|^2dx.
\end{equation}

The Cauchy problem for (\ref{rvp}) in the classical case ($c=+\infty$) is subcritical in the energy space and smooth initial data, say $f_0\in \calC_0^1$ (compactly supported $\calC^1$ functions), yield global in time  solutions, see Pfaffelmoser \cite{Pf}, Lions, Perthame \cite{LP}. The key to global existence is the uniform bound on the kinetic energy which follows from the interpolation estimate: 
\be
\label{gn2}
|\nabla\phi_f|_{L^2}^2\leq K\,||v|^2f|_{L^1}^{1/2}\,|f|_{L^1}^{\frac{7p-9}{6(p-1)}}\,|f|_{L^p}^{\frac{p}{3(p-1)}},
\end{equation}
for $9/7 < p < +\infty$. In the relativistic case $0<c<+\infty$, the Cauchy theory of smooth solutions is so far restricted to smooth radial data $$f_0\in {\cal{C}}^1_{0,rad}=\{f:\RR^{2N}\to [0,+\infty) \ \ \mbox{with radial symmetry, compactly supported and} \ {\cal{C}}^1\},$$ see Glassey, Schaeffer  \cite{GS}, Kiessling, Tahvildar-Zadeh \cite{T}, and then again a uniform bound on the kinetic energy suffices to ensure global existence. However, the relativistic problem is critical according to the interpolation estimate:
 \be
\label{gn}
|\nabla\phi_f|_{L^2}^2\leq K\,||v|f|_{L^1}\,|f|_{L^1}^{\frac{2p-3}{3(p-1)}}\,|f|_{L^p}^{\frac{p}{3(p-1)}},
\end{equation}
for $3/2<p<+\infty$. Glassey and Schaeffer proved in \cite{GS} that radially symmetric solutions to \fref{rvp} (for
 $c<+\infty$) with negative Hamiltonian blow up in finite time. In \cite{LMR4}, a {\it stable} self similar blow up dynamic for the relativistic problem corresponding to a concentration phenomenon is fully described.

We address in this paper the question of the existence and the stability of ground states for  the relativistic problem. In the classical case, this question has attracted considerable attention.  A large class of stationary solutions has been constructed in  \cite{Batt}, among which the functions of the microscopic energy $F(\frac{|v|^2}{2}+\phi)$. 
The question of stability of such steady states has been addressed in many works in the past and still stimulates a number of research programs. The first work on this subject goes back to Antonov (1960') \cite{A1,A2}, where a stability criterion for polytropes was established for the {\em linearized} GVP equation.
Then, this linear stability has been improved to a non linear stability by  Wolansky \cite{wol} for the so-called polytropes and later extended by  Guo \cite{Guo1} and, Guo and Rein \cite{GR1} to more general steady states. These analyses are based on the construction of steady states as minimizers of one-parameter energy-Casimir functionals:
\be
\label{energy-casimir}
\inf_{|f|_{L^1}=M} \calH_{\infty}(f) +\int_{\R^6} j(f) dx dv, 
\ee 
where $j$ is a suitable strictly convex function on $\R_+$. 
Extensions of this approach can also be found in \cite{Guo2,GR3}, completed
by Schaeffer \cite{S}. Non variational approaches based on linearization techniques have also been explored in \cite{Wan} and, more recently in \cite{GR4}.

In \cite{LMR-note,LMR1}, we observed that the frame of concentration compactness techniques as introduced by Lions \cite{PLL1}, \cite{PLL2} directly allows one to derive the existence of a {\it two parameters} family of ground states --in accordance with the scaling symmetry of the problem-- corresponding to the minimization problem:
\be
\label{2c}
I(M_1,M_j)=\inf_{|f|_{L^1}=M_1, \ |j(f)|_{L^1}=M_j} \calH (f), \ \ M_1,M_j>0
\end{equation} 
for a large class of convex functions $j$. In the case of polytropes $j(t)=t^p$, this has also been independently observed by S\'anchez and Soler \cite{SS}.

 
 \subsection{Statement of the results}
 
 
Here we propose to extend the variational approach of \cite{LMR-note,LMR1} to the relativistic framework. Very few papers have been devoted to stability issues in this setting. A stability result of some steady states solutions has been obtained recently by Had\v zi\'c and Rein \cite{Ha-Re} for specific perturbations. Our main claim in this paper is that following Lieb, Yau \cite{Lieb}, the breaking of the scaling symmetry in the relativistic case allows one to derive a similar variational theory of ground states like in the classical case under an additional subcritical size assumption. This will lead to a stability theory of ground states in the full energy space.

To wit, consider a strictly convex and even function $j:\R\to \R^+$ satisfying the following assumptions.

\noindent (B1)\  $j$ is a $\calC^1$ strictly convex function with  $j(0)=j'(0)=0.$
 
\noindent  (B2)\ There exists $p>3/2$ such that:
   
  \be 
  \label{j-sup-p}j(t)\geq C t^p, \ \ \forall t\geq 0. 
  \end{equation}
 
\noindent (B3)\ There exist $p_1, p_2 > 3/2$, such that
\begin{equation}
\label{j-dich-simple}
p_1\leq \frac{tj'(t)}{j(t)}\leq p_2\,\qquad \forall t>0.
\end{equation}
Note that this assumption (B3) is equivalent to the usual dichotomy condition
\be
\label{j-dichotomie}
b^{p_1}j(t) \leq j(bt) \leq b^{p_2}j(t), \ \ \forall \ b\geq 1, \ t\geq 0.
\end{equation}
Indeed, \fref{j-dichotomie} is equivalent to
$$\forall b\geq 1, \quad\forall t\geq 0, \qquad (bt)^{-p_1}j(bt) \geq t^{-p_1}j(t),\qquad (bt)^{-p_2}j(bt) \leq t^{-p_2}j(t),$$
which means that $t^{-p_1}j(t)$ is nondecreasing and $t^{-p_2}j(t)$ is nonincreasing on $\R_+$. Taking the derivative yields \fref{j-dich-simple}.

For a function $j$ satisfying (B1), (B2), (B3), we define the corresponding energy space
\be
\label{energyspace-j}
\calE_j =\{f\geq 0\ \ \mbox{with} \ \ |f|_{\calE_j}=|f|_{L^1}+|j(f)|_{L^1}+|\gamma_c(v)f|_{L^1}<+\infty\}.
\end{equation}
>From the interpolation inequality \fref{gn}, one can define a positive constant $K_j$ by
 \be
 \label{def-Kj}
 K_j= \inf _{f\in \calE_j-\{0\}} \frac{||v|f|_{L^1}\,|f|_{L^1}^{\frac{2p-3}{3(p-1)}}\,|j(f)|_{L^1}^{\frac{1}{3(p-1)}}}{|\nabla\phi_f|_{L^2}^2},
 \end{equation}
 where $\calE_j$ is the relativistic energy space \fref{energyspace-j}.

 \begin{proposition}[Existence of relativistic ground states]
 \label{th1}
 Let $j$ be a real function satisfying assumptions (B1), B2), (B3). Let $M_1>0, M_j>0$, and $c\in ]0,+\infty]$ be such that
\be
\label{cond-M1Mj}
M_1^{\frac{2p-3}{3(p-1)}}\,M_j^{\frac{1}{3(p-1)}}< 2c K_j,
\end{equation}
where $K_j$ is defined by \fref{def-Kj}. Let $\calF (M_1,M_j)= \{ f\in \calE_j, \ |f|_{L^1}=M_1, \ \ |j(f)|_{L^1}=M_j\}$. Then every minimizing sequence of the problem:
\be
\label{2c-rel}
I_c(M_1,M_j)= \inf_{f\in \calF (M_1,M_j)} \calH_c(f),
\end{equation}
where $\calH_c$ is defined by \fref{defhamiltonian},
is relatively strongly compact in the energy space $\calE_j$ up to a translation shift in the space variable $x$. Moreover, a  minimizer  $Q_j$ of \fref{2c-rel} is of the form 
\be
\label{forme-min} Q_j(x,v)= (j')^{-1} \left( \frac{\gamma_c(v)+\phi_{Q_j}(x) -\lambda}{\mu}\right)_+,\ \ \lambda,\mu<0,
\end{equation}
where $\phi_{Q_j}$ is linked to $Q_j$ according to \fref{phif}, has radial symmetry up to a translation shift and is a $\calC^2$ function on $\RR^3$. Eventually,  $Q_j$ is a compactly supported stationary solution to (\ref{rvp}).
\end{proposition}
Here we denote $t_+=\max(t,0)$ for all $t\in \RR$.
\begin{Rk} Note that the criterion (\ref{cond-M1Mj}) is optimal in the sense that for $M_1^{\frac{2p-3}{3(p-1)}}\,M_j^{\frac{1}{3(p-1)}}>2c K_j$, one can prove from scaling argument that $I_c(M_1,M_j)=-\infty$. Remark also from the definition \fref{def-Kj} of $K_j$ and from the conservation laws \fref{loi1} that an initial data $f_0\in \calC^1_{0,rad}$ satisfying
$$
|f_0|_{L^1}^{\frac{2p-3}{3(p-1)}}\,|j(f_0)|_{L^1}^{\frac{1}{3(p-1)}}< 2c K_j
$$
provides a global in time solution to \fref{rvp}.
\end{Rk}

Following Cazenave, Lions \cite{CL}, Proposition \ref{th1} classically implies the orbital stability in the energy space of {\it the family of ground states $\{Q_j\}$} only, and not of each $Q_j$ individually. Indeed a standard difficulty occurs here related to the uniqueness of the minimizer to (\ref{2c-rel}) --up to the symmetries--. In the classical case and for one-parameter variational problems, this issue was overcome in the analysis by Guo and Rein \cite{GR3} thanks to a uniqueness result due to Schaeffer \cite{S}. However, uniqueness of the minimizer, even only locally, in the general framework of Proposition \ref{th1} is a problem of independent importance and up to now mostly open.\\

We now claim that the orbitally stability of $Q_j$ may be derived even though uniqueness is not known using extra rigidities of the flow provided by the nonlinear transport. In particular, {\it equimeasurability properties} of the flow will allow us to prove some local isolatedness of the $Q_j$, and the following theorem follows which completes the analysis of orbital stability of \cite{LMR1} and extends it to the relativistic case:

\begin{theorem}[Orbital stability of classical and relativistic ground states]
\label{th2}
 Let $j$ be a real function satisfying assumptions (B1), B2), (B3). Let $M_1>0$, $M_j>0$ and $c\in ]0,+\infty]$ be such that \fref{cond-M1Mj} holds. Then any minimizer $Q_j$ of \fref{2c-rel}
is orbitally stable under the flow \fref{rvp} in  the energy space \fref{energyspace-j}. More precisely, 
given $\eps>0$, there exists $\delta(\eps)>0$ such that the
following holds true. 

i) Classical case $c=+\infty$. Let $f_0\in {\cal{C}}^1_0$
with $|f_0-Q_j|_{\calE_j}\leq \delta(\eps)$,  and let f(t) be the
classical solution to  (\ref{rvp}) with initial data $f_0$, then there exists a translation shift $x(t)\in\R^3$ such that $\forall t\in[0,+\infty)$, $$ \left|f(t,x+x(t),v)-Q_j\right|_{\calE_j}<\eps\,.$$

ii) Relativistic case $0<c<+\infty$. Let $f_0\in {\cal{C}}^1_{0,rad}$ 
with $|f_0-Q_j|_{\calE_j}\leq \delta(\eps)$,  then the classical solution f(t) to (\ref{rvp}) with initial data $f_0$ is global in time and $\forall t\in[0,+\infty)$, $$ \left|f(t)-Q_j\right|_{\calE_j}<\eps\,.$$
\end{theorem}
\begin{Rk}
In the classical case $c=+\infty$, the result of Theorem 1 can also be formulated in the framework of weak (or renormalized) solutions, as the Cauchy theory in this case is completely understood \cite{HH,DPL1,DPL2}.
\end{Rk}

Together with the results of \cite{LMR4} on blow up dynamics, Theorem \ref{th2} shows that the gravitational Vlasov Poisson system displays at least two stable dynamics: a global dynamic near ground state type solutions with a size strictly below the critical size required for blow up, a stable self similar blow up dynamic which may occur just above the critical threshold. This situation seems to be the generic one for critical problems with broken scaling invariance, see for example \cite{PR} for similar results in the context of nonlinear Schr\"odinger equations. Let us also mention also the works \cite{EL1},  \cite{EL2} on the pseudo-relativistic Boson star equation which are somehow connected both on the physical and mathematical side and where stable subcritical ground states are exhibited while a critical type finite time blow up regime is expected for larger masses. 

\section{Existence of relativistic ground states}


This section is devoted to the proof of Proposition \ref{th1} which is a consequence of the standard concentration compactness techniques, \cite{PLL1}. We shall adapt the proof of \cite{LMR1}.


\subsection{Properties of the infimum}


Let us start by summarizing some monotonicity properties
of the infimum \fref{2c-rel}.

 \begin{proposition}[Monotonicity properties of the infimum  $I_c(M_1,M_j)$]
 \label{prop1}
 Let $j$ be a real function satisfying assumptions (B1)-(B3) and let $M_1>0, M_j>0$ and $c\in ]0,+\infty]$ be such that \fref{cond-M1Mj} holds. Let $I_c(M_1,M_j)$ be the infimum defined by  
 \fref{2c-rel}, then we have:  
 \be
\label{IM1}
 -\infty < I_c(M_1,M_j) < 0
 \end{equation}
 and there holds the nondichotomy condition:  for all $0< \alpha <1$ and $0\leq \beta \leq 1$,
 \be
\label{I-dich}
I_c(\alpha M_1,\beta M_j)+ I_c((1-\alpha) M_1,(1-\beta) M_j) > I_c(M_1,M_j).
\end{equation}
 \end{proposition}
 
\begin{proof}
The case $c=+\infty$ is treated in \cite{LMR1,LMR-note} and we therefore assume $0<c<+\infty$.\\

{\bf Step 1.} The infimum is negative.\\

We first prove \fref{IM1}.  Let $\calH_c$ given by (\ref{defhamiltonian}) and $f\in \calE_j$ with
$|f|_{L^1}= M_1$ and $|j(f)|_{L^1}= M_j$, then from the definition \fref{def-Kj} of $K_j$, we have:
\be
\label{lowerboundH}
\begin{array}{ll} \ds \calH_c(f) &\ds \geq \left(\int_{\R^6} \gamma_c(v)
  f\right) -  \frac{1}{2K_j} \left(\int |v| f\right)
  M_1^{\frac{2p-3}{3(p-1)}}\,M_j^{\frac{1}{3(p-1)}}\\
&\ds \geq  \left(\int_{\R^6} \gamma_c(v)
  f\right) -  \frac{1}{2 c K_j} \left(\int (\gamma_c(v) +c^2) f\right)
  M_1^{\frac{2p-3}{3(p-1)}}\,M_j^{\frac{1}{3(p-1)}}\\
&\ds \geq \left(\int_{\R^6} \gamma_c(v)
  f\right) \left(1-
  \frac{1}{2 c K_j}M_1^{\frac{2p-3}{3(p-1)}}\,M_j^{\frac{1}{3(p-1)}}\right)
  -\frac{c}{2K_j}M_1^{\frac{5p-6}{3(p-1)}}\,M_j^{\frac{1}{3(p-1)}} \\
&\ds \geq -\frac{c}{2K_j}M_1^{\frac{5p-6}{3(p-1)}}\,M_j^{\frac{1}{3(p-1)}},
\end{array}
\end{equation}
where we have used condition \fref{cond-M1Mj} and the fact that
$\sqrt{1+|v|^2/c^2} \geq |v|/c$. This proves that  $I_c(M_1,kM_j)$ is bounded from below. To prove that it is negative, we use a rescaling argument. Assume that $f$ is moreover compactly supported and let $\tilde f(x,v)= f(\frac{x}{\lambda}, \lambda v)$, $\lambda >0$, then:

\begin{eqnarray*}
\calH_c (\tilde f) & =  & \frac{1}{\lambda} \left( \int_{\R^6}c^2\left(\sqrt{\lambda^2 +|v|^2/c^2} -\lambda\right) f - \frac{1}{2} |\nabla\phi_f|_{L^2}^2\right)\\
& = &  \frac{1}{\lambda}\int_{\R^6}\frac{|v|^2f}{\sqrt{\lambda^2 +\frac{|v|^2}{c^2}} + \lambda} - \frac{1}{2\lambda} |\nabla\phi_f|_{L^2}^2\\
& \sim & - \frac{1}{2\lambda} |\nabla\phi_f|_{L^2}^2 \ \ \mbox{as} \ \ \lambda\to +\infty
\end{eqnarray*}
and  \fref{IM1} follows.\\

{\bf Step 2.}  Monotonicity properties of the infimum.\\

We now claim the following monotonicity properties: for all $0<k\leq 1$,
 \be
 \label{mon-Mj}
 I_c(M_1,kM_j) \geq k^{\frac{1}{3(p_2-1)}} I_c(M_1,M_j),
 \end{equation}
 \be
 \label{mon-M1} I_c(kM_1,M_j) \geq k^{\frac{5p_1-6}{3(p_1-1)}} I_c(M_1,M_j).
 \end{equation}
Proof of \fref{mon-Mj}. We fix a  real number $0<k\leq 1$  and consider $f\in \calE_j$ such that
$|f|_{L^1}= M_1$ and $|j(f)|_{L^1}= kM_j$.  We introduce the rescaled function
$$\tilde f(x,v) = \alpha f(\alpha^{1/3} x, v)
,$$
then $|\tilde f|_{L^1}=M_1$ and $|j(\tilde f)|_{L^1}=h(\alpha,f)kM_j$ where:
\be
\label{defhalpha}
h(\alpha,f)= \frac{|j(\alpha f)|_{L^1}}{\alpha |j(f)|_{L^1}}.
\end{equation} 
Observe that $h(1,f)=1$ and $h(\alpha,f)\to +\infty$ as $\alpha\to +\infty$ from \fref{j-dichotomie}, hence from $k\leq 1$, we can find $\alpha\geq 1$ such that
$$h(\alpha,f)= \frac{1}{k} \ \ \mbox{and thus} \ \ |j(\tilde f)|_{L^1}=M_j.$$ Moreover, from $\alpha\geq 1$ and \fref{j-dichotomie}:
\be 
\label{borne-gamma1}
\frac{1}{k}\leq \alpha^{p_2-1}\,.
\end{equation}
We now compute
$$\calH_c (\tilde f)= \int_{\R^6} \gamma_c(v) f - \frac{1}{2} \alpha ^{1/3}|\nabla\phi_f|_{L^2}^2,$$ and thus from \fref{borne-gamma1}:
$$I_c(M_1,M_j)\leq \calH_c (\tilde f)\leq \int_{\R^6} \gamma_c(v) f - \frac{1}{2} \left(\frac{1}{k}\right)^{\frac{1}{3(p_2-1)}}|\nabla\phi_f|_{L^2}^2 \leq  \left(\frac{1}{k}\right)^{\frac{1}{3(p_2-1)}} \calH_c (f),$$
where we have used that $k\leq 1$. This concludes the proof of \fref{mon-Mj}.\\
Proof of \fref{mon-M1}. Similarly, we take $f$ such that $|f|_{L^1}= kM_1$ and $|j(f)|_{L^1}= M_j$ and set $$\tilde f(x,v) = \alpha f(\alpha^{1/3} k^{1/3} x, v).$$
We have $|\tilde f|_{L^1}= M_1$ and $|j(\tilde f)|_{L^1}=h(\alpha,f)k^{-1}M_j$. Now from (B1), $h(0,f)=0$ and we may find $\alpha\leq 1$ such that
$$h(\alpha,f)=k \ \ \mbox{and thus} \ \ |j(\tilde f)|_{L^1}=M_j.$$ Moreover, from $\alpha \leq1$ and \fref{j-dichotomie}:
$$
k\leq \alpha^{p_1-1}
$$
from which:
$$\begin{array}{ll}\ds I_c(M_1,M_j)\leq \calH_c (\tilde f)& \ds = 
k^{-1} \int_{\R^6}  \gamma_c(v)f - \frac{1}{2}k^{-5/3} \alpha^{1/3}|\nabla\phi_f|_{L^2}^2\\
&\ds \leq 
k^{-1} \int_{\R^6} \gamma_c(v) f - \frac{1}{2}k^{-5/3} k^{\frac{1}{3(p_1-1)}}|\nabla\phi_f|_{L^2}^2\\
&\ds \leq  k^{-\frac{5p_1-6}{3(p_1-1)}} \calH_c (f),
\end{array}$$
where we have used the fact that $\frac{5p_1-6}{3(p_1-1)} >1$ from $p_1>3/2$.  Hence \fref{mon-M1} follows.\\

{\bf Step 3.} The nondichotomy property.\\

Let $0<\alpha<1$, $0\leq \beta\leq 1$, then \fref{mon-Mj} and \fref{mon-M1} imply:
$$I_c(\alpha M_1,\beta M_j) \geq 
\alpha^{\frac{5p_1-6}{3(p_1-1)}} \beta^{\frac{1}{3(p_2-1)}} I_c( M_1, M_j),$$
and a similar inequality by exchanging $\alpha$ and $\beta$
with $1-\alpha$ and $1-\beta$ respectively.
As $I_c( M_1, M_j) <0$, \fref{I-dich} is equivalent to
the following
$$ \alpha^{\frac{5p_1-6}{3(p_1-1)}} \beta^{\frac{1}{3(p_2-1)}} + 
(1-\alpha)^{\frac{5p_1-6}{3(p_1-1)}} (1-\beta)^{\frac{1}{3(p_2-1)}}< 1$$
 which holds true since
 $ p_2>1$ and  $\frac{5p_1-6}{3(p_1-1)} >1 $ (since $p_1> 3/2$). This concludes the proof of  Proposition \ref{prop1}.
\end{proof}


\subsection{The concentration compactness argument}

We now prove Proposition \ref{th1}. We adapt the argument from \cite{LMR1, LMR-note} to which we refer for more details.\\

{\noindent{\bf Proof of Proposition \ref{th1}}

\bs
\ni
{\bf Step 1.} Compactness of the minimizing sequences.\\

 Let $M_1,M_j>0$, then from Lemma \ref{prop1} we know that $I_c(M_1, M_j)$ is finite and negative. Take then a minimizing sequence $f_n$ of \fref{2c-rel}:
\be
\label{suite-min}
|f_n|_{L^1}= M_1, \ \  |j(f_n)|_{L^1}= M_j \ \  \mbox{and} \ \ \lim_{n\to +\infty} \calH(f_n)= 
I_c(M_1, M_j).
\end{equation}
Let
$$\rho_n(x)= \int_{\R^3} f(x,v)dv.$$
We know from the concentration compactness principle developed in 
\cite{PLL1,PLL2}  that  there exists a subsequence $\rho_{n_k}$ for
which  one of the three possibilities occurs ($B_R$ being the ball of radius $R$ centered at the origin in $\R^3$): 

\begin{itemize}
\item {\em Compactness}: there exists $y_k\in \R^3$ such that 
\be
\label{concentration}
\forall \e>0, \ \ \exists R<+\infty \ \ \mbox{such that for all} \ \ k\geq 1\quad \int_{y_k+B_R}\rho_{n_k}(x)dx\geq M_1-\e\,;
\end{equation}
\item {\em Vanishing}: 
$$
\forall  R<+\infty, \ \ \lim_{k\to+\infty}\sup_{y\in\R^3}\int_{y+B_{R}}\rho_{n_k}(x)dx=0 \,;
$$
\item{\em Dichotomy:} there exists $\alphadef\in(0,M_1)$ such that for all
$\e>0$, there exist  subsequences
$(\rho_k^1)_{k\geq 1}, \
(\rho_k^2)_{k\geq 1}\in  L_+^1(\R^{3})$ and $k_0\geq 1$ such that for all $k\geq k_0$, 
$$ \ \  \left   \{ \begin{array}{llll}
\rho_{n_k}=\rho_k^1+\rho_k^2+w_k \ \ \mbox{with} \ \ 0\leq \rho_k^1,\rho_k^2,w_k\leq \rho_{n_k}, \ \ \rho_k^1\rho_k^2=\rho_k^1w_k=\rho_k^2w_k=0 \ \ \mbox{a.e.},\\[2mm]
dist(Supp(\rho_k^1),Supp(\rho_k^2))\to +\infty \ \ \mbox{as} \ \ k\to+\infty,\\[2mm]
\left|\rho_{n_k}-\rho^1_k-\rho^2_k\right|_{L^1}\leq \e, \ \ \left|\int_{\R^3} \rho^1_k(x)dx-\alphadef\right|+\left|\int_{\R^3} \rho^2_k(x)dx-(M_1-\alphadef)\right|<\e.
\end{array}
\right .
$$
\end{itemize}
We claim that only compactness can occur. Indeed, if vanishing occurs, then from Lemma 3.1 in \cite{LMR1}, 
$$|\nabla\phi_{f_{n_k}}|_{L^2}\rightarrow 0, \ \ as \ \ k \rightarrow 0.$$
Passing to the limit into
$$\calH_c (f_{n_k})= \int_{\R^6}  \gamma_c(v) f_{n_k}- \frac{1}{2} 
|\nabla\phi_{f_{n_k}}|_{L^2}^2 \geq - |\nabla\phi_{f_{n_k}}|_{L^2}^2,$$
leads to $I_c(M_1, M_j) \geq 0$ which contradicts property \fref{IM1}. Dichotomy cannot occur as it would violate the nondichotomy property \fref{I-dich}, see again \cite{LMR1} for further details.

We conclude that the compactness occurs on a subsequence. Now observe from \fref{suite-min}, (B2)  that $f_n$ is a bounded sequence in $L^p \cap L^1$. Moreover,  the kinetic energy $\int_{\R^6}  \gamma_c(v)f_n(x,v) dxdv$ is also bounded  thanks to the lower bound \fref{lowerboundH}. This implies the compactness in the velocity variable $v$. We deduce that we have $L^1$ compactness in $x$ up to a translation shift and that no concentration can occur in $x,v$ due to the $L^p$ boundedness. Thus the Dunford-Pettis criterion ensures:
 $$f_{n_k} (.+y_k) \rightharpoonup f \ \ \mbox{in} \ \ L^1,L^p$$
and from (\ref{concentration}):
 $$ \int_{\R^6} f(x,v) dx dv = M_1.$$
 This implies from a standard compactness argument --see \cite{LMR1}--:
  $$|\nabla\phi_{f_{n_k}}|_{L^2}^2 \rightarrow |\nabla\phi_{f}|_{L^2}^2,$$ and thus by lower semi-continuity of the $L^p$ norms:
\be 
\label{lsc}
\calH_c(f) \leq I_c( M_1, M_j), \ \  \ \ |j(f)|_{L^1}\leq  M_j.
\end{equation}
Then  from \fref{mon-Mj},

$$I_c(M_1,M_j)\geq \calH_c(f)\geq I_c( M_1, |j(f)|_{L^1})\geq \left(\frac{|j(f)|_{L^1}}{M_j}\right)^{\frac{1}{3(p_2-1)} }I_c(M_1, M_j) .
$$
This implies from \fref{IM1} that $|j(f)|_{L^1}\geq M_j$ and thus, from \fref{lsc}, we obtain
$$|f|_{L^1}= M_1 \ \ \  \mbox{and} \ \ \ |j(f)|_{L^1}= M_j$$ 
which together with \fref{lsc} implies that $f$ is a minimizer of 
\fref{2c-rel}. Moreover, we get: $$|f_{n_k}|_{L^1}\rightarrow |f|_{L^1},\ \ \ |\gamma_c(v)(f_{n_k})|_{L^1} \rightarrow |\gamma_c(v) f|_{L^1},\ \ \ 
|j(f_{n_k})|_{L^1} \rightarrow |j(f)|_{L^1}. $$ 
We now conclude from standard convexity arguments, see \cite{LMR1}, \cite{Brezis-Lieb}, that $f_{n_k} (.+y_k) \to f$ in $L^1$, $|v|^2f_{n_k} (.+y_k) \to |v|^2f $ in $L^1$ and $j(f_{n_k} (.+y_k) -f))\to 0$ in $L^1$, hence the strong convergence in the energy space $\calE_j$.\\

{\bf Step 2.} Euler-Lagrange equation for the minimizer.\\

Let $Q_j$ be a minimizer of \fref{2c-rel}. Let $Q_j^*$ be the nondecreasing symmetric rearrangement of $Q_j$ in $x$ --see \cite{Lieb-Loss} for instance--, then $|Q_j^*|_{L^1}=M_1$, $|j(Q_j^*)|_{L^1}=M_j$ and $\calH_c(Q_j^*) \leq \calH_c (Q_j)$,  this inequality being strict unless $\phi_{Q_j}$ is radial in $x$ up to a space translation shift. From now on, without loss of generality, we assume that $\phi_{Q_j}$ is radial around $0$. From standard Euler-Lagrange theory --see \cite{LMR1} for further details--, there exist  constants $\lambda$ and $\mu$ such that:
\be
\label{euler-lagrange1}
\gamma_c(v) + \phi_{Q_j}(x) = \lambda + \mu j'(Q_j), \ \ \ \mbox{on \ the support \ of} \ Q_j\,.
\end{equation}

We now claim that $\lambda,\mu<0$. We first multiply \fref{euler-lagrange1} by $Q_j$ and integrate over $(x,v)\in\R^6$:
\be
\label{mf}
\int_{\R^6} \gamma_c(v) Q_j + \int_{\R^3} \phi_{Q_j} Q_j = \lambda M_1+ \mu \int_{\R^6} j'(Q_j)Q_j.
\end{equation}
Next, we multiply \fref{euler-lagrange1} by $v\cdot \nabla_vQ_j$ and integrate by parts to get:
\be
\label{mnablavf}
\int_{\R^6} \gamma_c(v) Q_j + \int_{\R^3} \phi_{Q_j} Q_j  +\frac{1}{3} 
\int_{\R^6} \frac{|v|^2}{\sqrt{1+|v|^2/c^2}} Q_j= \lambda M_1+ \mu M_j.
\end{equation}
Similarly, we multiply \fref{euler-lagrange1} by $x\cdot \nabla_xQ_j$ and get:
\be
\label{mnablaxf}
\int_{\R^6} \gamma_c(v) Q_j + \frac{5}{6}\int_{\R^3} \phi_{Q_j} Q_j  = \lambda M_1+ \mu M_j.
\end{equation}
Combining \fref{mnablavf} with \fref{mnablaxf}, we get the relativistic Viriel identity:
\be
\label{viriel}
\int_{\R^6} \frac{|v|^2}{\sqrt{1+|v|^2/c^2}} Q_j= -\frac{1}{2}\int_{\R^3} \phi_{Q_j} Q_j.
\end{equation}
Therefore, we have
\be
\label{Inegatif}
I_c(M_1,M_j) = \calH_c(Q_j)=  -c^2\int_{\R^6} \left( 1- \frac{1}{\sqrt{1+|v|^2/c^2}}\right) Q_j.
\end{equation}
Let
\be
E_{kin}(Q_j)= \int_{\R^6} \gamma_c(v) Q_j,
\end{equation}
then, we have from \fref{viriel}
\be
\label{viriel2}
\int_{\R^6} \frac{|v|^2}{\sqrt{1+|v|^2/c^2}} Q_j= E_{kin}(Q_j) -I_c(M_1,M_j).
\end{equation}
Now, reporting \fref{viriel2} and \fref{viriel} in relations \fref{mf} and 
\fref{mnablaxf}, we get
\be
\label{EL1}
- E_{kin}(Q_j) + 2 I_c(M_1,M_j)= \lambda M_1+ \mu \int_{\R^6} j'(Q_j)Q_j
\end{equation}
and
\be
\label{EL2}
- \frac{2}{3} E_{kin}(Q_j) +\frac{5}{3}  I_c(M_1,M_j)= \lambda M_1+ \mu M_j 
\end{equation}
Subtracting  these  two last  relations gives:

\be
\label{muj}
3 \mu \int_{\R^6} \left( j'(Q_j)Q_j -j(Q_j)\right)= I_c(M_1,M_j)- E_{kin}(Q_j) < 0
\end{equation}
and thus the strict convexity of $j$ implies $\mu <0$.
Now reporting the expression of $\mu$ in \fref{EL2} leads to
\be
\label{lambda}
\begin{array}{l}\ds 3 \lambda M_1 \int_{\R^6} \left(j'(Q_j)Q_j -j(Q_j)\right)= \\
\ds 
\ \ \ \ - E_{kin}(Q_j)\int_{\R^6} \left( 2j'(Q_j)Q_j -3j(Q_j)\right) +I_c(M_1,M_j)
\int_{\R^6} \left( 5j'(Q_j)Q_j -6j(Q_j)\right),
\end{array}
\end{equation}
thus $\lambda<0$ from (B3).\\

{\bf Step 3.} Regularity of the potential and compact support of $Q_j$.\\

Let us now prove the $\calC^2$ regularity of $\phi_{Q_j}$ on $[0,+\infty)$. Using the Euler-Lagrange equation \fref{euler-lagrange1}, we have
$$\rho_{Q_j}(x)= \int_{\R^3} (j')^{-1} \left( \frac{\gamma_c(v) + \phi_{Q_j}(x) -\lambda}{\mu}\right)_+ dv.$$
We then pass to the spherical velocity coordinate $u=|v|$ and perform the change of variable $q=\gamma_c(v)/|\mu|$, to get
\begin{equation}
\label{h1}
\rho_{Q_j}(x)=4\pi\int_0^{+\infty} (j')^{-1} \left(\frac{\phi_{Q_j}(r)-\lambda}{\mu}-q\right)_+ |\mu|c \left(1+\frac{|\mu| q}{c^2}\right)
\left[\left(1+\frac{|\mu| q}{c^2}\right)^2-1\right]^{1/2} dq.
\end{equation}
Using this expression, we shall now control $\rho_{Q_j}$ by a power of the potential
$\phi_{Q_j}$. We first recall that the support of $Q_j$ is contained in the set of $(x,v)$ such that
\be
\label{suppQ_j}
\gamma_c(v) + \phi_{Q_j}(r) -\lambda \leq 0,  \ \ \mbox{with} \ \ r=|x|.
\end{equation}
This in particular implies that $0\leq \lambda - \phi_{Q_j}(r) \leq |\phi_{Q_j}(r)|$ and 
$|\phi_{Q_j}| \geq |\lambda|>0$, 
on the support of $\rho_{Q_j}$.
>From \fref{j-dich-simple} we also  have $(j')^{-1}(s)\leq Cs^{1/(p-1)}$,
hence we straightforwardly get
from \fref{h1}
$$\begin{array}{ll}
\ds \rho_{Q_j}(x)&\ds \leq K \int_0^{+\infty} \left(\frac{\phi_{Q_j}(r)-\lambda}{\mu}-q\right)_+^{\frac{1}{p-1}}(1+q^2) dq.\\[5mm]
&\ds \leq K\left(\frac{\phi_{Q_j}(r)-\lambda}{\mu}\right)_+^{1+\frac{1}{p-1}} +K\left(\frac{\phi_{Q_j}(r)-\lambda}{\mu}\right)_+^{3+\frac{1}{p-1}}\leq K|\phi_{Q_j}(r)|^{3+\frac{1}{p-1}}
\end{array}
$$
on the  support of $\rho_{Q_j}$ (recall indeed that $|\phi_{Q_j}|\geq
|\lambda|>0$ on this support). Now we follow a bootstrap argument as
in \cite{LMR4} (proof of Proposition 2): if $\rho_{Q_j}\in L^{q_k}$,
then $\phi_{Q_j}=\frac{1}{4\pi|x|}\star \rho_{Q_j}\in L^{r_k}$ with
$r_k=\frac{3q_k}{3-2q_k}$ from Hardy, Littlewood, Sobolev so that
$\rho_{Q_j}\in L^{q_{k+1}}$ with
$$q_{k+1}\left(3+\frac{1}{p-1}\right)=r_k \ \ \mbox{ie} \ \
q_{k+1}=\frac{3(p-1)q_k}{(3p-2)(3-2q_k)}.$$ A simple analysis of the
sequence $q_k$ with $q_0=\frac{6}{5}$ shows that $p>\frac{3}{2}$
implies that there exists $k_0=k_0(p)$ such that $q_{k_0}>\frac{3}{2}$
and thus $\rho_{Q_j}\in L^{q_{k_0}}$. From Sobolev embeddings, this
implies that $\phi_{Q_j}\in \calC^{0,\alpha}$ for some
$0<\alpha<1$. Finally, by \fref{h1}, $\rho_{Q_j}$ is continuous on
$[0,+\infty)$, which is enough to deduce from  the Laplace equation
\begin{equation}
\label{lap}
(r^2\phi_{Q_j}'(r))'= r^2\rho_{Q_j}(r),
\end{equation}
 that $\phi_{Q_j}$ is a strictly increasing $\calC^2$ function on
 $[0,+\infty)$.  Furthermore, as $\rho_{Q_j}\in L^1\cap L^{q_{k_0}}$ with
 $q_{k_0}>\frac{3}{2}$, the potential
 $\phi_{Q_j}=\frac{1}{4\pi|x|}\star \rho_{Q_j}$ goes to $0$ at
 infinity.

To complete the proof, we now show that $Q_j$ is compactly
supported. Observe from \fref{suppQ_j} and the monotonicity of
$\phi_{Q_j}$ on $[0,+\infty[$ that $\gamma_c(v) \leq  \lambda -\phi_{Q_j}(0)$
and  $\phi_{Q_j}(r) \leq \lambda <0$ on the support of $Q_j$.  These two
last inequalities, together with the fact that $\phi_{Q_j}(r)$ goes to
$0$ as $r\rightarrow +\infty$,   imply that the support of $Q_j$ is compact in $(x,v)$.
\qed


\section{Orbital stability of the ground states}


We now prove the orbital stability of the ground states $Q_j$ without the knowledge of the uniqueness of the minimizer for the problem (\ref{2c-rel}). The key is a uniqueness statement of the minimizer under a condition of equimeasurability which is inherited from the transport evolution, see Proposition \ref{prop2}.


\subsection{Proof of theorem \ref{th2}}
 

We treat the case $c<\infty$ which forces us to restrict to compactly supported smooth radial solutions due to the Cauchy theory. The case $c=+\infty$ would be treated similarly without the radial assumption restriction.

\bs

{\noindent{\bf Proof of Theorem \ref{th2}}

\bs
\ni
Let $Q_0$ be a radial minimizer of \fref{2c-rel} and assume that the result of Theorem \ref{th2} is false. Then there exist $\varepsilon >0$ and sequences $f_0^n\in \calC^1_{0,rad}$\,, $t_n>0$, 
such that 
\be
\label{f0nQn}
\lim_{n\to +\infty}  \left|f_0^n-Q_0\right|_{\calE_j} =0,
\end{equation}
 and
\be
\label{absurde1}
\forall n\geq 0, \ \ \left|f^n(t_n,x,v)-Q_0\right|_{\calE_j}\geq \eps,
\end{equation}
where $f^n(t,x,v)$ is a solution to \fref{rvp} with initial data
$f_0^n$.
This means in particular that
\be
\label{f0n}
\lim \calH_c (f_0^n)=I_c(M_1,M_j) ,\qquad \lim |f_0^n|_{L^1}= M_1,
\qquad \lim |j(f_0^n-Q_0)|_{L^1}=0.
\end{equation}
In particular, $f_0^n$ converges to $Q_0$ in the
strong $L^p$ topology and hence
almost everywhere, up to a subsequence.  From assumptions (B1)--(B3) and the convexity of $j$, this implies form classical argument
 --see Theorem 2 in Br\'ezis and Lieb \cite{Brezis-Lieb}--:
\be
\label{saturationj}
|j(f_0^n)|_{L^1}-|j(Q_0)|_{L^1} \rightarrow 0\ \ \mbox{as} \ \
n\rightarrow +\infty.
\end{equation}
Let now $g_n(x,v)=f^n(t_n,x,v)$. Then from the conservation properties of the Vlasov-Poisson flow,  there holds $$\lim_{n\to +\infty} {\calH_c}(g_n)= I_c(M_1,M_j) , \ \ |g_n|_{L^1}=
M_1, \qquad |j(g_n)|_{L^1}= M_j\,,$$
which means that $g_n$ is a minimizing sequence of \fref{2c-rel}. From Proposition \ref{th1},  $g_n$ is relatively strongly compact in $\calE_j$: 
\be
\label{gnconverge}
g_n\to Q_1 \ \ \mbox{in} \ \ \calE_j
\end{equation}
for some minimizer $Q_1$ --without any translation due to the radial assumption--. Moreover, for  all smooth compactly supported $\theta$, there holds the conservation law:
\be
\label{egalite-n}\int_{\RR^6} \theta (g_n) = \int_{\RR^6} \theta (f_0^n),
\end{equation}
and hence passing to the limit $n\to +\infty$ from \fref{f0nQn} and \fref{gnconverge}, we get:
\be
\label{hiehi}
\int_{\RR^6} \theta (Q_1) = \int_{\RR^6} \theta (Q_0).
\end{equation}
Now, from standard arguments, \fref{hiehi} implies the equimeasurability of $Q_1$ and $Q_0$\,:
$$\forall t>0, \ \ meas\{Q_0(x,v)>t\}=meas\{Q_1(x,v)>t\}.$$

Let now $h_n=f^n(t'_n,x,v)$ where $t_n'$ is the time such that 
\begin{equation}
\label{e2}
\left|f^n(t'_n,x,v)-Q_0\right|_{\calE_j}=\frac{\eps}{2},
\end{equation}
which is well defined by the continuity of the flow. Then arguing as above for $Q_1$, $h_n\to Q_2$ in $\calE_j$ where $Q_2$ is a minimizer of \fref{2c-rel} which satisfies the equimeasurability property: 
\be
\label{equinoe}
\forall t>0, \ \ meas\{Q_0(x,v)>t\}=meas\{Q_1(x,v)>t\}=meas\{Q_2(x,v)>t\}.
\end{equation}
In particular, from \fref{absurde1}  and \fref{e2}, $\{Q_0,Q_1,Q_2\}$ are three distinct solutions to the Euler-Lagrange equation \fref{forme-min}  for some a priori distinct Euler-Lagrange multipliers $(\lambda_j, \mu_j)_{0\leq j\leq 2}$, and satisfying (\ref{equinoe}). This now contradicts the following uniqueness statement of equimeasurable steady solutions to \fref{rvp} which is the core of our argument.

\sloppy
\begin{proposition}[Isolatedness of equimeasurable steady states]
\label{prop2}
Let $c\in]0,+\infty]$  and $G$ be a given  continuous and strictly increasing function on $\R_+$ such that $G(0)=0$. Let $\calA$ be  the set of all functions of the form:
\be
\label{hiwuof}
Q(x,v) =G\left( \frac{\gamma_c(v) +\phi(x)
  -\lambda}{\mu}\right)_+,   \ \ \forall (x,v) \in \R^6,\
  \end{equation}
where $\lambda <0, \ \mu<0$ are arbitrary constants and where $\phi$
is a ${\cal{C}}^2$ radial solution to the elliptic equation:
\be
\label{laplace}
\Delta \phi = \int_{\R^3}G\left(\frac{\gamma_c(v) +\phi(x)
  -\lambda}{\mu}\right)_+ dv= \rho_Q(x), \ \ \phi(r)\to 0\ \ \mbox{as}\ \ r\to +\infty .
\end{equation}
Let  $Q_0 \in \calA$, then the following set
\be
\label{calS} \calS(Q_0) =\left\{ Q\in \calA,\  s.t. \ \ meas\left\{Q(x,v) >t\}=
   meas\{Q_0(x,v) >t\right\}, \ \ \forall t>0\right\},
\end{equation}
cannot contain more than two elements. 
\end{proposition}

Applying Proposition \ref{prop2} with $G=(j')^{-1}$ yields a contradiction and concludes the proof of Theorem \ref{th2}.
\qed


\subsection{Isolatedness of equimeasurable steady state solutions}
 

We now turn to the proof of Proposition \ref{prop2}. Note that the difficulty is that the Lagrange multiplier $(\lambda,\mu)$ in \fref {hiwuof} are a priori assumed to be different for the solutions we consider. If they were the same, uniqueness in the ODE sense for \fref{hiwuof} would immediately conclude the proof (see below). Even in the more restricted setting of solutions being minimizers of \fref{2c-rel}, the Lagrange multipliers cannot be simply connected to quantities conserved by the flow due to the breaking of the scaling symmetry --this would be the case for $c=+\infty$ where the proof can therefore be simplified--.

\bs
{\noindent{\bf Proof of Proposition \ref{prop2}}

\bs
\ni
Let $Q_0,Q$ be two radially symmetric solutions of \fref{hiwuof} for some parameters $(\lambda_0,\mu_0)<0$, $(\lambda,\mu)<0$.\\

{\bf Step 1.} Comparison of the $L^{\infty}$ norms.\\

 Let $\phi_0,\phi$ be the Poisson potentials associated respectively to $Q_0,Q$. We claim the first relation:
 \be
\label{egalite1}
\left(\frac{\phi(0)
  -\lambda}{\mu}\right)_+=\left( \frac{\phi_0(0)
  -\lambda_0}{\mu_0}\right)_+. 
\end{equation}
 Proof of  \fref{egalite1}:  As $\phi_0$ is radial, we shall use the abuse of notation
$\phi_0(x)=\phi_0(r)$ with $r=|x|$, and write the Laplace equation $\Delta
\phi_0 = \rho_0(x)= \int_{\R^3} Q_0(x,v)dv,$ as
$$r^2\phi_0'(r)=\int_0^rs^2\rho_0(s)ds.$$
In particular this implies that $\phi_0$ is a nondecreasing function in
$r$ and we have: $\phi_0(r) \geq \phi_0(0)$. Now, since $G$ is nondecreasing and since $\mu <0$, we have from \fref{hiwuof}:
\be
\label{Q0infini}
\forall (x,v) \in \R^6, \ \ Q_0(x,v) \leq G \left( \frac{\phi_0(0)
  -\lambda_0}{\mu_0}\right)= |Q_0|_{L^{\infty}}.
\end{equation}
We now observe that 
\be
\label{dhio}
|Q|_{L^{\infty}}=|Q_0|_{L^{\infty}}.
\end{equation}
Indeed, if $|Q_0|_{L^{\infty}}<|Q|_{L^{\infty}}$, we may choose $t$ such that $|Q_0|_{L^{\infty}}<t<|Q|_{L^{\infty}}$
to get
$$meas\{Q_0(x,v)>t\}=0 \ \ \mbox{and} \ \ meas\{Q(x,v)>t\}>0$$ contradicting \fref{calS}. \fref{Q0infini}) and (\ref{dhio}) now imply \fref{egalite1}.\\

{\bf Step 2.} Two possible values for $\mu$\\

>From \fref{egalite1} and (\ref{Q0infini}), \be
\label{defa}
a=\frac{\phi(0)
  -\lambda}{\mu}=\frac{\phi_0(0)
  -\lambda_0}{\mu_0}>0
\end{equation}
or otherwise $Q=Q_0=0$ and the proof is over. From \fref{calS}, $$meas\{Q_0(x,v)>t\}=meas\{Q(x,v)>t\}>0, \ \ \forall t \in [0,G(a)],$$
which we rewrite equivalently:
$$meas\left\{ \gamma_c(v) + \phi_0(r)-\phi_0(0) < |\mu_0|(a-\tau)\right\}= meas\left\{ \gamma_c(v) + \phi(r)-\phi(0) < |\mu|(a-\tau)\right\},
$$
for all $\tau \in [0,a]$. We now identify the leading terms in this equality
when $\tau \rightarrow a$. First remark that the set $\{\gamma_c(v) +
\phi_0(r)-\phi_0(0) < |\mu_0|(a-\tau)\}$
goes to $\{0\}$ when $\tau\rightarrow a$. Therefore,  by expanding 
$\gamma_c(v)$ near $v=0$ and $\phi_0(r)$ near $r=0$ and taking the
      leading terms, we get
$$
meas\left\{\gamma_c(v) + \phi_0(r)-\phi_0(0) < |\mu_0|(a-\tau)\right\} \sim
      meas\left\{\frac{|v|^2}{2}  + \phi_0''(0)\frac{|x|^2}{2} < |\mu_0|(a-\tau)
      \right\},
$$ when $t\rightarrow a$.
The measure (in $\R^6$)  of the above rhs can be computed
      explicitly leading to 
\be
\label{meas-equiv2} meas\left\{\gamma_c(v) + \phi_0(r)-\phi_0(0) < |\mu_0|(a-\tau) \right\} \sim
      K \left(\frac{-\mu_0}{\sqrt{\phi_0''(0)}}\right)^3 (a-\tau)^3,
\end{equation}
for some universal constant $K>0$, and thus:
\be
\label{egalite2}
\frac{\phi''(0)}{\mu^2}= \frac{\phi_0''(0)}{\mu_0^2}.
\end{equation}
Now, we shall use the elliptic equations satisfied by $\phi_0$ and
$\phi$. The Laplace equation \fref{laplace} for $\phi$ can be
written is spherical coordinates as 
$$ \phi''(r) + 2\frac{\phi'(r)}{r}=\int_{\R^3}G \left(\frac{\gamma_c(v) +\phi(r)
  -\lambda}{\mu}\right)_+ dv.$$
 Taking $r=0$ in this equation and using \fref{defa}, we get:
$$3 \phi''(0)=\int_{\R^3}G \left(\frac{\gamma_c(v)}{\mu} + a
 \right)_+ dv.$$ Let us now focus onto the case $c<+\infty$.
Passing to the spherical velocity coordinate $u=|v|$ and performing the
 change of variable 
$$q= \frac{\gamma_c(u)}{|\mu|},$$
we obtain:
$$3 \phi''(0)=4\pi\int_0^{+\infty}G \left((a-q)_+\right) |\mu|c \left(1+\frac{|\mu| q}{c^2}\right)
\left[\left(1+\frac{|\mu| q}{c^2}\right)^2-1\right]^{1/2} dq.$$
A similar identity holds for $\phi_0$ and $\mu_0$. Combining this with
\fref{egalite2}, we get that $\mu$ must be a solution to
\be
\label{egaliteF}
 F(|\mu|)= F(|\mu_0|),
\end{equation} where $F$ is defined on $]0,+\infty[$ by  
\be
\label{defF}
F(s)= c\int_0^{a}\frac{1}{s}
\left(1+\frac{s q}{c^2}\right)
\left[\left(1+\frac{s q}{c^2}\right)^2-1\right]^{1/2} G(a-q)dq.
\end{equation}
We compute the derivatives of $F(s)$ and find
\be
\label{Fseconde}
F''(s)= \int_0^{a} \frac{q^2}{c^3}\left( \frac{q}{c^2}+\frac{3}{s}\right)
\left[\left(1+\frac{s q}{c^2}\right)^2-1\right]^{-3/2} G(a-q)dq>0,
\end{equation}
and thus $F$ is strictly convex on $\R_+$.
Therefore, for given $\mu_0$, equation \fref{egaliteF} in $|\mu|$ has at most two solutions and thus $\mu_0,\mu<0$ ensures that $\mu$ can take at most two values.
\begin{Rk}
If $c=+\infty$,  then $F$ given by \fref{defF} simplifies into $F(s) = K s^{-1/2}$ for some universal constant $K>0$. In particular, $F$ is 
a nonincreasing function in this case and thus \fref{egaliteF} implies $\mu=\mu_0$. In this case, the set $\calS(Q_0)$ defined by \fref{calS} is reduced to $\{Q_0\}$.
\end{Rk}

{\bf Step 3.} Conclusion\\

We now claim from an ODE type argument that $\mu_0=\mu$ implies $Q_0=Q$. Indeed, let then
$$\psi_0(r)= \phi_0(r) - \lambda_0,  \ \ \mbox{and} \ \ \psi(r)= \phi(r) - \lambda,$$
then from  
\fref{defa} and the radial symmetry,  we have
\be
\label{cauchy-psi}
\psi(0) = \psi_0(0), \ \  \mbox{and} \ \ \psi'(0) = \psi_0'(0)=0. 
\end{equation}
Moreover, $\psi$ and $\psi_0$ solve the same radial Laplace equation (since $\mu=\mu_0$):
\begin{equation}
  \label{ode45}
\Delta \psi = \int_{\R^3}G \left(\frac{\gamma_c(v) +\psi(r)}{\mu}\right)_+ dv=h(\psi(r)),
\end{equation}
with
$$h(u)=4\pi\int_0^u \left(1+\frac{u-q}{c^2}\right)\left[\left(1+\frac{u-q}{c^2}\right)^2-1\right]^{1/2}G\left(\frac{q}{\mu}\right)dq.$$
It is clear that $h$ is a $\calC^1$ function of $u$ thus, together with \fref{cauchy-psi}, uniqueness for the ODE \fref{ode45} implies that $\psi=\psi_0$. Hence $Q=Q_0$ from (\ref{hiwuof}). This concludes the proof of Proposition \ref{prop2}.
\qed

\begin{acknowledgement} P. Rapha\"el was supported by the Agence Nationale de la Recherche, ANR ONDENONLIN. M. Lemou was supported by the Agence Nationale de la Recherche, ANR Jeunes Chercheurs MNEC. F. M\'ehats was supported by the Agence Nationale de la Recherche, ANR project QUATRAIN. 
\end{acknowledgement}


\begin{thebibliography}{10}

\bibitem{A1}
Antonov, A. V., Remarks on the problem of stability in stellar dynamics. {\em Soviet Astr., AJ.,} {\bf 4},
859-867 (1961).

\bibitem{A2} Antonov, A. V., Solution of the problem of stability of a stellar system with the Emden density law and spherical velocity distribution. {\em J. Leningrad Univ.Se. Mekh. Astro. } {\bf 7}, 135-146 (1962).


%
\bibitem{Batt} Batt, J.; Faltenbacher, W.; Horst, E., Stationary spherically symmetric models in stellar dynamics, Arch. Rat. Mech. Anal. 93, 159-183 (1986).

%
\bibitem{binney} Binney, J.; Tremaine, S., Galactic Dynamics, Princeton University Press, 1987.

%

%

%
\bibitem{Brezis-Lieb} Br\'ezis, H.;  Lieb, E.: A relation between pointwise convergence of functions and convergence of functionals. Proceedings of the American Mathematical Society. Volume 88, no 3, July 1983. 

%
\bibitem{CL}  Cazenave, T.; Lions, P.-L. Orbital stability of standing waves for some nonlinear Schr\"odinger equations, Comm. Math. Phys. 85 (1982), no. 4, 549--561.

%


%
\bibitem{DPL1} DiPerna, R. J.; Lions, P.-L., Solution globale de type Vlasov-Poisson.  C. R. Acad. Sci. Paris Sér I math. 307 (1988), no 12, 655-658. 

%
\bibitem{DPL2} DiPerna, R. J.; Lions, P.-L., Global weak solutions of kinetic equations.  Rend. Sem. Mat. Univ. Politec. Torino  46  (1988),  no. 3, 259--288 (1990).


%
 

%

%
\bibitem{fridman} Fridmann, A. M.; Polyachenko, V. L., Physics of gravitating systems, Springer-Verlag

%
\bibitem{EL1}  Fr\"ohlich, J.; Jonsson, B.; Lars G.; Lenzmann, E., Boson stars as solitary waves, Comm. Math. Phys. 274 (2007), no. 1, 1--30.

\bibitem{EL2} Fr\"ohlich, J.; Lenzmann, E., Blowup for nonlinear wave equations describing boson stars, Comm. Pure Appl. Math. 60 (2007), no. 11, 1691--1705


%
\bibitem{GS} Glassey, R.T.; Schaeffer, J., On symmetric solutions of the relativistic Vlasov-Poisson system, Comm. Math. Phys. 101 (1985), no. 4, 459--473.

%




\bibitem{Guo1} Guo, Y., Variational method for stable polytropic galaxies, Arch. Rat. Mech. Anal.  130 (1995), 163-182. 

\bibitem{Guo2} Guo, Y., On the generalized Antonov's stability criterion.
{\em Contem. Math.}  {\bf 263}, 85-107 (2000)

\bibitem{GR1} Guo, Y.; Rein, G., Stable steady states in stellar dynamcics, Arch. Rat. Mech. Anal.  147 (1999), 225--243.




\bibitem{GR3} Guo, Y.; Rein, G., Isotropic steady states in galactic dynamics, Comm. Math. Phys. 219 (2001), 607--629.

\bibitem{GR4} Guo, Y.; Rein, G., A non-variational approach to nonlinear stability in stellar dynamics applied to the King model.  Comm. Math. Phys.  271  (2007),  no. 2, 489--509.
   
\bibitem{Ha-Re} Had\v zi\'c, M.; Rein, G., Global existence and nonlinear stability for the relativistic Vlasov-Poisson system in the gravitational case , 
Indiana Univ. Math. J. 56, 2453-2488 (2007).

\bibitem{HH} Horst, E.; Hunze, R., Weak solutions of the initial value problem for the unmodified nonlinear Vlasov equation.
Math. Methods Appl. Sci. 6. (1984), no. 2, 262-279.

\bibitem{T} Kiessling M. K. H, Tahvildar-Zadeh A. S, On the relativistic Vlasov-Poisson system,  arXiv preprint 0708.1760, to appear in Indiana Univ. Math. J.

%

%
\bibitem{LMR1} Lemou, M.; M\'ehats, F.; Rapha\"el, P., On the orbital
  stability of the ground states and the singularity formation for the
  gravitational Vlasov-Poisson system, to appear in Arch. Rational
  Mech. Anal.

\bibitem{LMR-note} Lemou, M.; M\'ehats, F.; Rapha\"el, P.,  Orbital stability and singularity formation for Vlasov-Poisson systems.  C. R. Math. Acad. Sci. Paris  341  (2005),  no. 4, 269--274. 


%
\bibitem{LMR2} Lemou, M.; M\'ehats, F.; Rapha\"el, P., Structure of the linearized gravitational Vlasov-Poisson system close to a polytropic ground state, to appear in SIAM J. Math. Anal.

%

\bibitem{LMR4} Lemou, M.; M\'ehats, F.; Rapha\"el, P., Stable self-similar blow-up dynamics for the three dimensional  gravitational Vlasov-Poisson system,
To appear in J. Amer. Math. Soc. 

%
\bibitem{Lieb-Loss}Lieb, E. H.; Loss, M., Analysis. Second edition. Graduate Studies in Mathematics, 14. American Mathematical Society, Providence, RI, 2001.
 
%
\bibitem{Lieb} Lieb, E.H.; Yau, H.T., The Chandrasekhar theory of stellar collapse as the limit of quantum mechanics, Comm. Math. Phys. 112 (1987), no. 1, 147--174.

\bibitem{LP} Lions, P.-L.; Perthame, B., Propagation of moments and regularity for the $3$-dimensional Vlasov-Poisson system.  Invent. Math. 105 (1991), no. 2.

%
\bibitem{PLL1} Lions, P.-L., The concentration-compactness principle in the calculus of variations. The locally compact case. I. Ann. Inst. H. Poincaré Anal. Non Linéaire 1 (1984), no. 2, 109--145. 

%
\bibitem{PLL2} Lions, P.-L., The concentration-compactness principle in the calculus of variations. The locally compact case. II. Ann. Inst. H. Poincaré Anal. Non Linéaire 1 (1984), no. 4, 223--283.




%
%
%
%

%
\bibitem{Pf} Pfaffelmoser, K.; Global classical solutions of the
  Vlasov-Poisson system in three dimensions for general initial data,
  J. Diff. Eq. {\bf 95} (1992), 281-303.

 %
 \bibitem{PR} Planchon, F.; Rapha\"el, P., Existence and stability of the log-log blow-up dynamics for the $L\sp 2$-critical nonlinear Schr\"odinger equation in a domain, Ann. Henri Poincar\'e 8 (2007), no. 6, 1177--1219

%
%
%


\bibitem{SS} S\'anchez, \'O.; Soler, J., Orbital stability for polytropic galaxies.  Ann. Inst. H. Poincaré Anal. Non Linéaire  23  (2006),  no. 6, 781--802.


%
\bibitem{S} Schaeffer, J., Steady States in Galactic Dynamics, Arch. Rational, Mech. Anal. {\bf 172} (2004), 1--19.

%
\bibitem{VKF} Van Kampen, N.G.; Felderhof, B.V., Theoretical methods in plasma physics, Amsterdam, North Holland 1967.

\bibitem{Wan} Wan, Y-H, On nonlinear stability of isotropic models in stellar dynamics, Arch. Ration. Mech. Anal 147 (1999), no.3, 245-268.
%

%
\bibitem{wol}  Wolansky, G., On nonlinear stability of polytropic galaxies. {\em Ann. Inst. Henri Poincaré,} {\bf 16}, 15-48 (1999).


\end{thebibliography}
\end{document}